\documentclass[12pt]{article}
\usepackage[cp1251]{inputenc}
\usepackage[russian]{babel}
\usepackage[pdftex, colorlinks, urlcolor=blue, pdfborder={0 0 0 [2 0]}]{hyperref}

\usepackage{amsfonts,amsmath,amssymb,theorem}
\usepackage{array}
\usepackage{esint}
\usepackage{graphicx}

\newtheorem{theorem}{Теорема}

\newtheorem{remark}{Замечание}
\newtheorem{dfn}{Определение}

\newenvironment{proofs}
{\vspace{1pt}\par{\sl%
Доказательство.\,\ }}%
{\noindent\vspace{1pt}}

\renewcommand{\refname}{Библиографический список}

\title{О двух подходах к получению обобщенной формулы Ито -- Вентцеля\footnote{%
\textcolor{blue}{\textit{Препринт  № 174 / В.~А.~Дубко, Е.~В.~Карачанская ; Вычисл. центр ДВО РАН.~-- Хабаровск : Изд-во Тихоокеан. гос. ун-та, 2012; ISBN 978~-~5~-~7389~-~1059~-~3.}}}}

\author{Дубко В. А.\footnote{%
Авиационный национальный университет, г. Киев, Украина}, Карачанская Е. В.\footnote{%
Тихоокеанский госуниверситет, г. Хабаровск, Россия; {\it e-mail}: EKarachanskaya@mail.khstu.ru}
}
\begin{document}
\date{}

\maketitle

\textit{AMS 2010 классификация} Primary: 60G60, Secondary: 60H15; 60H30

\textit{Ключевые слова}: стохастическое дифференциальное уравнение Ито со скачками, формула Ито -- Вентцеля, ядро интегрального инварианта, аппроксимация $\delta-$функции.

\begin{abstract}
Рассмотрено два подхода к построению обобщенной формулы Ито -- Вентцеля:  первый -- при непосредственном привлечении обобщенной формулы Ито,  второй~-- на основе представления о ядрах интегральных инвариантов.
\end{abstract}

\tableofcontents

\section*{Предварительные замечания}
\addcontentsline{toc}{section}{Предварительные замечания \dotfill}

Пусть ${\bf x}(t)\in \mathbb{R}^{n} $ решение уравнения
\begin{equation*}
\frac{d{\bf x}(t)}{dt} =a(t),\ \ \ \ {\bf x}(t)|_{t=0} ={\bf x}(0).
\end{equation*}
Пусть функция $F(t;{\bf x})$ непрерывна вместе со своими частными производными $F'_{t} (t;{\bf x})$, $F'_{x_{i} } (t;{\bf x})$, $i=\overline{1,n}$. Тогда, опираясь на  правило дифференцирования сложной функции, получаем:
\begin{equation*}
\frac{dF(t;{\bf x}(t))}{dt} =F'_{t} (t;{\bf x}(t))+a_{j} (t)F'_{x_{j} } (t;{\bf x}(t)),\ \ \ \ F(t;{\bf x})|_{t=0} =F(0;{\bf x}).
\end{equation*}
Если же $F(t;{\bf x})$ есть решение уравнения
\begin{equation*}
\frac{\partial F(t;{\bf x})}{\partial t} =g(t;{\bf x}), \ \ \ \ \ F(t;{\bf x})|_{t=0} =F(0;{\bf x}) ,
\end{equation*}
то переходим к выражению:
\begin{equation*}
\frac{dF(t;{\bf x}(t))}{dt} =g(t;{\bf x}(t))+a_{j} (t)F'_{x_{j} } (t;{\bf x}(t)).
\end{equation*}
Рассмотрим случай, когда ${\bf x}(t)$ определяется уравнением
\begin{equation} \label{GrindEQ__1_1_2_}
 d{\bf x}(t)=a(t;{\bf x})dt+B_{k} (t;{\bf x})dw_{k} (t),
\end{equation}
а $F(t;{\bf x})$ подчинено уравнению:
\begin{equation} \label{GrindEQ__1_4_1_} dF(t;{\bf x})=Q(t;{\bf x})dt+D_{k} (t;{\bf x})dw_{k} (t),
\end{equation}
где ${\bf w}(t)$ -- тот же  $m$-мерный винеровский процесс с  независимыми  компонентами, как и в уравнении \eqref{GrindEQ__1_1_2_}.

Относительно случайных функций $a(t;{\bf x})$, $B_{k}(t;{\bf x})$, $Q(t;{\bf x})$ и $D_{k}(t;{\bf x})$ предполагаем, что они  непрерывны и ограниченные вместе со своими первыми  и вторыми частными производными по компонентам ${\bf x}$, неупреждающие по $t$ относительно приращений винеровского векторного процесса ${\bf w}(t)$ и согласованные с потоком $\sigma$-алгебр, индуцированным ${\bf w}(t)$ на всем интервале $[0,T]$. Этих ограничений достаточно для построения  дифференциала от процесса $F(t;{\bf x}(t))$ (см. теорему 1):
\begin{equation} \label{GrindEQ__1_4_2_}
\begin{array}{c}
d_{t} F(t;{\bf x}(t))=Q(t;{\bf x}(t))dt+D_k(t;{\bf x}(t))d{\rm w}(t)+ \\ + b_{i,k} (t)\displaystyle\frac{\partial }{\partial x_{i} } D_k(t;{\bf x}(t))\, dt + \\
+\displaystyle\left(a_{i} (t)\frac{\partial }{\partial x_{i} } +\frac{1}{2} b_{i,k} \left(t\right)b_{j,k} (t)\frac{\partial ^{2} }{\partial x_{i} \partial x_{j} } \right)   F(t;{\bf x}(t))dt+ \\
 +\displaystyle b_{i,k} (t)\frac{\partial }{\partial x_{i} } F(t;{\bf x}(t))dw_{k} (t), \end{array}
 \end{equation}
\begin{equation*}
F(t;{\bf x}(t))\left|_{t=0} =F(0;y),\right. \; F(0;y)\in C_{0}^{2} .
\end{equation*}
Соотношение  \eqref{GrindEQ__1_4_2_} носит название \textit{формулы Ито -- Вентцеля} \cite{Rozovsky_73}.

\begin{remark} В дальнейшем, как частный случай, формула \eqref{GrindEQ__1_4_2_} будет получена, когда присутствуют кроме винеровских и пуассоновские возмущения, с использованием  обобщенной формулы Ито, а также на основе уравнений для  ядер интегральных инвариантов.
\end{remark}
Цель данной работы в  построении подобной формулы  для случая,  когда случайный процесс ${\bf x}(t)$ является решением уравнения
\begin{equation*}
d{\bf x}(t)=a(t)dt+B(t)d{\bf w}(t)+\int\limits_{{\mathbb R}(\gamma)}g(t;\gamma) {\rm \; }\nu (dt;d\gamma ),
\end{equation*}
где ${\bf x}(t) \in \mathbb{R}^{n} $; ${\bf w}(t)$ -- $m$-мерный винеровский процесс; $\nu (dt;d\gamma )$ -- нецентрированная мера Пуассона,
и, соответственно,
\begin{equation*}
dF(t;{\bf x})=Q(t;{\bf x})dt+D_{k} (t;{\bf x})dw_{k} (t)+\int\limits_{{\mathbb R}(\gamma)} G(t;x;\gamma ) \nu (dt;d\gamma ).
\end{equation*}

Следуя традиции, полученную в дальнейшем формулу, будем называть \textit{обобщенной формулой Ито -- Вентцеля}.

Эту формулу получим как «прямым» методом -- на основе классической теории стохастических дифференциальных уравнений, так и опираясь на свойства ядер интегральных инвариантов.

Первоначально эта формула была получена именно на основе уравнений для ядер интегральных инвариантов,  в связи с поиском сохраняющихся функционалов для открытых систем \cite{D_02,11_KchOboz}.  В дальнейшем попытки получить эту формулу "<прямым"> методом предпринимались другими авторами, хотя и для частного случая, но, по нашему мнению, наиболее удачно в работе \cite{Oks_07}.
Именно последняя работа и стимулировала рассмотрение общего случая получения этой формулы.

Отметим, что в отличии  от известных работ других авторов, мы с самого начала будем работать c обобщенными уравнениями для нецентрированных мер. Это потребует несколько более жестких ограничений на класс весовых функций при пуассоновских возмущениях, но делает доказательство более компактным и прозрачным.

Для построения интегралов по центрированной пуассоновской мере достаточно, чтобы с вероятностью единица выполнялось условие \cite[c.\,254]{GS_68}:
\begin{equation*}
\int\limits_{{\mathbb R}(\gamma)}|f(\gamma )| ^{2} \Pi (d\gamma )<\infty.
\end{equation*}

При работе с нецентрированными мерами добавляется еще условие \linebreak \cite[c.\,255]{GS_68}:
\begin{equation*}
\int\limits_{{\mathbb R}(\gamma)}|f(\gamma )| \Pi (d\gamma )<\infty
\end{equation*}

\begin{remark} В дальнейшем, с целью компактности записи, вместо  $\int\limits_{{\mathbb R}(\gamma)} $ будем пользоваться обозначением $\int  $.
\end{remark}

\section*{1. Прямой метод получения обобщенной \\ формулы  Ито -- Вентцеля }
\addcontentsline{toc}{section}{1. Прямой метод получения обобщенной формулы \\ Ито -- Вентцеля  \dotfill}

Используются следующие  обозначения (\cite{GS_68}).

$H_{n}$ -- означает пространство случайных векторов (отображений) $\mathfrak{F}_t$- измеримых при произвольном $t\in [0,T]$ и с вероятностью 1  и таких, что (\cite[c.\,255]{GS_68}):
\begin{equation*}
\int\limits_{0}^{T}|\alpha (t)|^{n} dt<\infty.
\end{equation*}

Если  функция $\varphi (t;y)=\varphi (t;y;\omega )$, при произвольных $t\in [0,T]$ и $y\in \mathbb{R}(y)$, измеримая как функция трех аргументов $t;y;\omega $, и
\begin{equation*}
\int\limits_{0}^{T}\int  |\varphi (t;y)|^{2} \Pi (dy)dt<\infty  \end{equation*}
то используется обозначение: $\varphi (t;y)\in H_{2} (\Pi)$ \cite[c.\,257]{GS_68}.

Обобщенным дифференциалом Ито (обобщенным уравнением Ито) называют такое равенство (\cite[\S 6,\,c.\,256]{GS_68}):
\begin{equation}\label{GrindEQ__2_1_1_}
d{\bf x}(t)=a(t)dt+B(t)d{\bf w}(t)+\int g(t;\gamma )\nu (dt;d\gamma),
\end{equation}
где ${\bf x}(t)$, $a(t)=\{ a_{1} (t),...,a_{n} (t)\}$, $g(t;\gamma )=\{ g_{1} (t;\gamma ),...,g_{n} (t;\gamma )\} \in \mathbb{R}^{n} $, \linebreak ${\bf w}(t)$ -- $m$-мерный винеровский процесс;
\begin{equation} \label{GrindEQ__2_1_2_}
\begin{array}{c}
  |g(t;\gamma )|\in H_{1,2} (\Pi ),\ \ \ \ \sqrt{|a_{j} (t)|},\, |b_{j,k} (t)|\in H_{2} , \\
  B(t)=\{ {\rm \; }b_{j,k} (t); \ \ \ j=\overline{{1,n}};\ \ \ k=\overline{{1,m}},\}
\end{array}
\end{equation}
$\nu (\Delta t;\Delta \gamma )$ - стандартная пуассоновская мера на $\mathbb{R}^{n'} \times [0,T]$, моделирующая независимые случайные величины на интервалах и множествах, которые не пересекаются в пространстве $\gamma \in {\mathbb R}(\gamma)=\mathbb{R}^{n'} $.

\begin{theorem} \label{t1}
Пусть $F(t;{\bf x})$, $(t;{\bf x})\in[0;T]\times\mathbb{R}^{n} $ -- скалярная функция, дифференциал которой имеет вид:
\begin{equation} \label{GrindEQ__2_5_1_}
\partial _{t} F(t;{\bf x})=Q(t;{\bf x})dt+D_{k} (t;{\bf x})dw_{k} (t)+\int G(t;{\bf x};\gamma ) \nu (dt;d\gamma )
\end{equation}
и для коэффициентов \eqref{GrindEQ__2_5_1_} выполнены условия:
$\mathcal{L.}a)$: $Q(t;{\bf x})$,$D_{k} (t;{\bf x})$, \linebreak $G(t;{\bf x};\gamma )\in\mathbb{R}$ в общем случайные функции измеримые относительно того же потока $\sigma$-алгебр $\Bigl\{\mathcal{F}\Bigr\}_{0}^{T}$,
$\mathcal{F}_{t_{1}}\subset\mathcal{F}_{t_{2}}$, $t_{1}<t_{2}$,, что и процессы $w(t)$, $\nu (t;\mathcal{A})$ для любого множества  $\mathcal{A}\in \mathfrak{B}$  из фиксированной борелевской $\sigma $-алгебры  {\rm(\cite[с. 266]{GS_68})};
и
$\mathcal{L.}b)$: c вероятностью единица  функции $Q(t;{\bf x}),D_{k} (t;{\bf x}),G(t;{\bf x};u)$ непрерывны и ограничены по совокупности переменных вместе со своими вторыми частными производными по компонентам вектора ${\bf x}$. Тогда, если случайный процесс ${\bf x}(t)$  подчинен системе  \eqref{GrindEQ__2_1_1_} и  выполняются ограничения \eqref{GrindEQ__2_1_2_}, то существует стохастический дифференциал:
\begin{equation} \label{GrindEQ__2_5_2_}
\begin{array}{c}
 \displaystyle d_{t} F(t;{\bf x}(t))=Q(t;{\bf x}(t))dt+D_{k} (t;{\bf x}(t))dw_{k} +\\
 +b_{i,k} (t)\displaystyle\frac{\partial }{\partial x_{i} } F(t;{\bf x}(t))dw_{k} + \\
+\displaystyle \Bigl[a_{i} (t)\frac{\partial }{\partial x_{i} } F(t;{\bf x}(t))+\frac{1}{2} b_{i,k} (t)b_{j,k} (t)\frac{\partial ^{{\kern 1pt} 2} F(t;{\bf x}(t))}{\partial x_{i} \partial x_{j} }\Bigr. + \\
+\Bigl.\displaystyle
b_{i,k} (t)\frac{\partial }{\partial x_{i} } D_{k} (t;{\bf x}(t))\Bigr]dt+ \\
+\displaystyle\int \Bigl[(F(t;{\bf x}(t)+g(t;\gamma ))-F(t;{\bf x}(t))\Bigr]\nu (dt;d\gamma ) +\\
+ \displaystyle\int G(t;{\bf x}(t)+g(t;\gamma )\nu (dt;d\gamma ). \end{array}
\end{equation}
\end{theorem}

Выражение \eqref{GrindEQ__2_5_2_} -- \textit{обобщеннная формулы Ито -- Вентцеля} на случай присутствия как винеровских, так и пуассоновских возмущений (для обобщенного СДУ Ито).

\begin{proofs} Доказательство базируется на определении $\delta $- функции как границы некоторой специальной функции. Пример использования дельта функции для построения интегральных инвариантов был дан в работе \cite{D_84}.

Определяют  $\delta $-функцию через  асимптотические свойства по разному.  Для наших целей  достаточно ограничиться рассмотрением случаем такого асимптотического представления \cite[c.\,48-49]{Gelfand_00}:
\begin{equation} \label{GrindEQ__2_5_3_}
\begin{array}{l}
f(t;x)=\displaystyle\int\limits_{-\infty }^{\infty }f(t;y) \delta (y-x)dy=
 \\
 =\displaystyle\lim\limits_{\varepsilon \downarrow 0} \frac{1}{\varepsilon \sqrt{2\pi } } \int\limits_{-\infty }^{\infty }f(t;y) \exp \left\{ -\frac{(y-x)^{2} }{2\varepsilon ^{2} } \right\} dy.
\end{array}
\end{equation}
Положим, что для функции $f(t;x)$ выполняется условие Липшица:
\begin{equation*}|f(t;y_{1} )-f(t;y_{2} )|\le L(t)|y_{1} -y_{2} |^{\varsigma } ,0<\varsigma \le 1
\end{equation*}
Поскольку $t$ в доказательствах рассматривается как фиксированное, то для упрощения в записи формул, параметр $t$ будем опускать. Т.е., вместо $f(t;x)$ будем использовать $f(x)$, а  вместо $L(t)$ -- обозначение  $L$.

Выполним замену переменных
\begin{equation*}(y-x)\varepsilon ^{-1} =z,{\rm \; \; \; }y=\varepsilon z+x .
\end{equation*}
Находим:
\begin{equation*}
\begin{array}{c}
  \displaystyle\frac{1}{\sqrt{2\pi } } \int\limits_{-\infty }^{\infty }f(\varepsilon z+x) \exp \left\{ -\frac{z^{2} }{2} \right\} dz=f(x)\frac{1}{\sqrt{2\pi } } \int\limits_{-\infty }^{\infty } \exp \left\{ -\frac{z^{2} }{2} \right\} dz+ \\
  \\
  +\displaystyle\frac{1}{\sqrt{2\pi } } \int\limits_{-\infty }^{\infty }\left(f(\varepsilon z+x) -f(x)\right)\exp \left\{ -\frac{z^{2} }{2} \right\} dz
\end{array}
\end{equation*}
Оценим последний интеграл:
\begin{equation*}
\begin{array}{c}
  \displaystyle\left|\frac{1}{\sqrt{2\pi } } \int\limits_{-\infty }^{\infty }\left(f(\varepsilon z+x) -f(x)\right)\exp \left\{ -\frac{z^{2} }{2} \right\} dz\right|\le \\
 \le \displaystyle\frac{L}{\sqrt{2\pi } } \int\limits_{-\infty }^{\infty }|(\varepsilon z+x) -x|^{\varsigma } \exp \left\{ -\frac{z^{2} }{2} \right\} dz\le
\\
\le \displaystyle\varepsilon L\frac{2}{\sqrt{2\pi } } \left[\int\limits_{0}^{1}z^{\varsigma } \exp \left\{ -\frac{z^{2} }{2} \right\} dz -\int\limits_{1}^{\infty }\exp \left\{ -\frac{z^{2} }{2} \right\}d\left( -\frac{z^{2} }{2}\right)\right]\le
\\
\le \displaystyle\varepsilon L\frac{2}{\sqrt{2\pi } } \Bigl[z^{\varsigma } |_{z=1} +1\Bigr]\le \varepsilon L\frac{4}{\sqrt{2\pi } }  .
\end{array}
\end{equation*}
Так что при $\varepsilon \to 0$, для произвольных функций $f(t;x)$, которые для любых $t\in [0,T]$ ограничены с вероятностью 1 и   удовлетворяют условию Липшица по компоненте  $x$, справедливо равенство \eqref{GrindEQ__2_5_3_}.

При доказательстве теоремы в выражениях под знаками интегралов будут присутствовать производные  от $\delta $-функции. Покажем, что и в этом случае справедливы оценки, подобные выполненным выше, на основе соответствующих предельных представлений. Продифференцируем \eqref{GrindEQ__2_5_3_}, применим интегрирование по частям:
\begin{equation*}
\begin{array}{c}
\displaystyle\frac{\partial }{\partial x} f(x)=\int\limits_{-\infty }^{\infty }f(y) \frac{\partial }{\partial x} \delta (y-x)dy=
 \\
  =\displaystyle\lim\limits_{\varepsilon \downarrow 0} \frac{1}{\varepsilon \sqrt{2\pi } } \int\limits_{-\infty }^{\infty }f(y) \frac{\partial }{\partial x} \exp \left\{ -\frac{(y-x)^{2} }{2\varepsilon ^{2} }\right\} dy=
  \\
  =-\displaystyle\lim\limits_{\varepsilon \downarrow 0} \frac{1}{\varepsilon \sqrt{2\pi } } \int\limits_{-\infty }^{\infty }f(y) \frac{\partial }{\partial y} \exp \left\{ -\frac{(y-x)^{2} }{2\varepsilon ^{2} }\right\} dy=\\
 =-\displaystyle\lim\limits_{\varepsilon \downarrow 0} f(y)\frac{1}{\varepsilon \sqrt{2\pi } } \exp \left\{ -\frac{(y-x)^{2} }{2\varepsilon ^{2} } \right\} \Bigl. \Bigr|_{-\infty }^{+\infty } +\\
+\displaystyle\lim\limits_{\varepsilon \downarrow 0} \frac{1}{\varepsilon \sqrt{2\pi } } \int\limits_{-\infty }^{\infty } \exp \left\{ -\frac{(y-x)^{2} }{2\varepsilon ^{2} } \right\} \frac{\partial }{\partial y} f(y)dy,\ \ \ \varepsilon >0.
\end{array}
\end{equation*}

При условии, что производная $f'_{y} (y)$ удовлетворяет условию Липшица,  для асимптотического представления производной от $\delta (x)$ доказательство совпадает с предшествующим. Аналогично устанавливается справедливость соотношений, где присутствуют вторые производные $\delta (x)$, и выполнено условие Липшица  для второй производной $f''_{y} (y)$. Обобщая использовавшиеся ограничения видим, что в дальнейшем необходимо потребовать  существование  и непрерывность $f'_{y} (y)$ и липшицевость для вторых производных.

Перейдем к доказательству соотношения \eqref{GrindEQ__2_5_2_} теоремы, воспользовавшись установленным выше свойствами  $\delta (x)$ на классе всех непрерывных и ограниченных с вероятностью единица функций:
\begin{equation} \label{GrindEQ__2_5_4_}
F(t,{\bf x}(t))=\displaystyle\int\limits_{\mathbb{R}^{n} }\prod_{i=1}^{n}\delta (y_{i} -x_{i} (t))  F(t,{\bf y})d\Gamma ({\bf y}),
\end{equation}
где $d\Gamma ({\bf y})$ -- элемент фазового объема \cite{D_89,D_02}.

Условия на коэффициенты оказываются достаточными для существования границы сумм, которые аппроксимируют в  среднеквадратичном соответствующие интегралы вида \eqref{GrindEQ__2_5_4_} при допредельном представлении  $\delta $-функций.

Пусть  дифференциал для $F(t,{\bf y})$ определяется \eqref{GrindEQ__2_5_1_}, а для  ${\bf x}(t)$ -- уравнением \eqref{GrindEQ__2_1_1_}. Формально применим обобщенную формулу Ито к подынтегральному выражению \eqref{GrindEQ__2_5_4_}:
$$
dF(t,{\bf y})\displaystyle\prod _{i=1}^{n}\delta (y_{i} -x_{i} (t)) =F(t,{\bf y})\left[a_{i} (t)\frac{\partial }{\partial x_{i} } \prod _{i=1}^{n}\delta \left(y_{i} -x_{i} (t)\right) \right.+
$$
$$
 +\left.\displaystyle\frac{1}{2} b_{i,k} (t)b_{j,k} (t)\frac{\partial ^{2} }{\partial x_{i} \partial x_{j} } \prod _{i=1}^{n}\delta (y_{i} -x_{i} (t)) \right]dt+
$$
$$
 +b_{i,k} (t)D_{k} (t;{\bf x})\displaystyle\frac{\partial }{\partial x_{i} } \prod _{i=1}^{n}\delta (y_{i} -x_{i} (t)) dt+
$$
$$
  +F(t,y)b_{i,k} \displaystyle\frac{\partial }{\partial x_{i} } \prod _{i=1}^{n}\delta (y_{i} -x_{i} (t))dw_{k} (t)+
  $$
  $$
 +\displaystyle\prod _{i=1}^{n}\delta (y_{i} -x_{i} (t)) \Bigl[Q(t;{\bf x})dt+D_{k} (t;{\bf x})dw_{k} (t)\Bigr]+
$$
$$
 +\displaystyle\int \Bigl[\prod _{i=1}^{n}\delta (y_{i} -x_{i} (t)-g_{i} (t,\gamma ))  (F(t;{\bf x})+G(t;\gamma ))\Bigr.-
$$
$$
 -\displaystyle\Bigl.\prod _{i=1}^{n}\delta (y_{i} -x_{i} (t)) F(t;{\bf x})\Bigr]\nu (dt,d\gamma ).
$$

Перейдем в правой части равенства от частных производных по компонентам вектора ${\bf x}$ к частным производных по компонентам  вектора ${\bf y}$. Тогда приходим к такому представлению выражения выше:
\begin{equation*}
\begin{array}{c}
dF(t,{\bf y})\displaystyle\prod _{i=1}^{n}\delta (y_{i} -x_{i} (t)) =F(t,{\bf y})\Bigl[-a_{i} (t)\frac{\partial }{\partial y_{i} } \prod _{i=1}^{n}\delta (y_{i} -x_{i} (t))\Bigr. +
 \\
 +\displaystyle\frac{1}{2} b_{i,k} (t)b_{j,k} (t)\frac{\partial ^{2} }{\partial y_{i} \partial y_{j} } \prod _{i=1}^{n}\delta (y_{i} -x_{i} (t)) \Bigl.\Bigr]dt-
 \end{array}
 \end{equation*}
 $$
 -b_{i,k} (t)D_{k} (t,{\bf y})\displaystyle\frac{\partial }{\partial y_{i} } \prod _{i=1}^{n}\delta (y_{i} -x_{i} (t)) dt-
 F(t;{\bf y})b_{i,k} \displaystyle\frac{\partial }{\partial y_{i} } \prod _{i=1}^{n}\delta (y_{i} -x_{i} (t)) dw_{k} (t)+
$$
$$
  +\displaystyle\prod _{i=1}^{n}\delta (y_{i} -x_{i} (t)) D_{k} (t,y)dw_{k} (t)+
 \displaystyle\prod _{i=1}^{n}\delta (y_{i} -x_{i} (t)) Q(t,{\bf y})dt+
 $$
\begin{equation*}
\begin{array}{c}
 +\displaystyle\int \Bigl[\Bigr. \prod _{i=1}^{n}\delta (y_{i} -x_{i} (t)-g(t,\gamma )) (F(t,{\bf y})+G(t;\gamma ))-
 \\
 -\displaystyle\prod _{i=1}^{n}\delta (y_{i} -x_{i} (t)) F(t,{\bf y})\Bigl.\Bigr]\nu (dt,d\gamma ).
\end{array}
 \end{equation*}

Возвратимся к интегрированию по пространству $\mathbb{R}^{n} ({\bf y})$, с учетом равенства \eqref{GrindEQ__2_5_4_}. Воспользовавшись операцией интегрирования по частям, с учетом допредельных свойств $\delta $-функций, находим:
\begin{equation*}
\begin{array}{c}
d_{t} F(t,{\bf x}(t))=\displaystyle\int\limits_{\mathbb{R}(y)}\Bigl\{\Bigr. \prod _{i=1}^{n}\delta (y_{i} -x_{i} (t))\Bigl [\Bigr. a(t)+b_{i,k} (t)D_{k} (t,{\bf y})F(t,{\bf y})+
 \\
 +\displaystyle\frac{1}{2} b_{i,k} (t)b_{j,k} (t)\frac{\partial ^{2} }{\partial y_{i} \partial y_{j} } F(t,y)\Bigl.\Bigr]dt+b_{i,k} (t)\frac{\partial }{\partial y_{i} } F(t,{\bf y})dw_{k} (t)+
\end{array}
\end{equation*}
$$
  +Q(t,{\bf y})dt+D_{k} (t,y)dw_{k} (t)+
  $$
\begin{equation*}
\begin{array}{c}
 +\displaystyle\int \Bigl[\Bigr. \prod _{i=1}^{n}\delta (y_{i} -x_{i} (t)-g_{i} (t,\gamma )) (F(t,{\bf y})+G(t;\gamma ))- \\
 -\displaystyle\prod _{i=1}^{n}\delta (y_{i} -x_{i} (t)) F(t,y)\Bigl.\Bigr]\nu (dt,d\gamma )\Bigl.\Bigr\} d\Gamma ({\bf y})
\end{array}
\end{equation*}
Учитывая правило интегрирования при наличии $\delta $-функций, приходим к равенству \eqref{GrindEQ__2_5_2_}.
\end{proofs}

\section*{2. Метод получения обобщенной формулы \\ Ито -- Вентцеля  на основе представления \\об  интегральных инвариантах}
\addcontentsline{toc}{section}{2. Метод получения обобщенной формулы Ито -- Вентцеля \\ на основе представления об интегральных инвариантах \dotfill}
Перейдем теперь к получению соотношения \eqref{GrindEQ__2_5_2_} на основе представления об интегральных инвариантах и уравнений для них.

Пусть ${\bf x}(t)$ -- динамический процесс, определенный на
$\mathbb{R}^{n}$, являющийся решением системы стохастических
дифференциальных уравнений
\begin{equation}\label{Ayd01.1}
\begin{array}{c}
  dx_{i}(t)= \displaystyle a_{i}\Bigl(t;{\bf x}(t) \Bigr)\, dt
  +\sum\limits_{k=1}^{m}
b_{i,k}\Bigl(t;{\bf x}(t)\Bigr)\, dw_{k}(t) +
\int g_{i}\Bigl(t;{\bf x}(t);\gamma\Bigr)\nu(dt;d\gamma), \\
  {\bf x}(t)={\bf x}\Bigl(t;{\bf x}(0) \Bigr) \Bigr|_{t=0}={\bf
  x}(0), \ \ \ \ \ \ i=\overline{1,n}, \ \ \ \ \ \ t\geq 0,
\end{array}
\end{equation}
где ${\bf w} (t)$ -- $m$-мерный винеровский процесс,
$\nu(t;\Delta\gamma)$ -- однородная по $t$ нецентрированная мера
Пуассона.

Относительно коэффициентов ${ a}(t;{\bf x})$, ${ b}(t;{\bf x})$  и
${g}(t;{\bf x};\gamma)$ будем предполагать, что они ограничены и
непрерывны вместе со своими производными
$
\displaystyle \frac{\partial a_{i}(t;{\bf x})}{\partial x_{j}}$,\,
$\displaystyle \frac{\partial b_{i,k}(t;{\bf x})}{\partial
x_{j}}$,\,
$\displaystyle \frac{\partial^{2} b_{i,k}(t;{\bf x})}{\partial
x_{j}\partial
x_{l}}$,\,$\displaystyle \frac{\partial g_{i}(t;{\bf
x};\gamma)}{\partial x_{j}}$,\, $i,j,l=\overline{1,n}
$
по совокупности переменных $(t;{\bf x};\gamma)$ -- ограничения $\mathfrak{L})$.

\subsection*{2.1. Стохастическое ядро стохастического \\ интегрального  инварианта}
\addcontentsline{toc}{subsection}{2.1. Стохастическое ядро стохастического интегрального инварианта}

Рассмотрим систему обобщенных стохастических дифференциальных \linebreak
уравнений Ито вида  (\ref{Ayd01.1}).
Пусть $\rho(t;{\bf x};\omega)$ -- случайная функция, измеримая
относительно потока $\sigma$-алгебр
$\Bigl\{\mathcal{F}\Bigr\}_{0}^{T}$,
$\mathcal{F}_{t_{1}}\subset\mathcal{F}_{t_{2}}$, $t_{1}<t_{2}$,
согласованного с процессами ${\bf w} (t)$ и $\nu(t;\Delta\gamma)$
(далее параметр $\omega$ будем опускать) и  относительно любой функции $f(t;{\bf x})$ из класса локально ограниченных функций, имеющей ограниченные вторые
производные по ${\bf x}$ для нее выполнены
соотношения:
\begin{equation}\label{Aydyad1}
\displaystyle\int\limits_{\mathbb{R}^{n}}\rho(t;{\bf x})f(t;{\bf
x})d\Gamma({\bf x})=\int\limits_{\mathbb{R}^{n}}\rho(0;{\bf
y})f(t;{\bf x}(t;{\bf y}))d\Gamma({\bf y})
\end{equation}
\begin{equation}\label{Aydyad2}
\displaystyle\int\limits_{\mathbb{R}^{n}}\rho(0;{\bf
x})d\Gamma({\bf x})=1,
\end{equation}
\begin{equation}\label{Aydyad3}
\begin{array}{c}
\displaystyle  \lim\limits_{|{\bf x}|\to \infty}\rho(0;{\bf x})=0,
\ \ \ d\Gamma({\bf x})=\prod\limits_{i=1}^{n}dx_{i},
\end{array}
\end{equation}
где ${\bf x}(t;{\bf y})$ -- решение уравнения (\ref{Ayd01.1}).

Если $f(t;{\bf x})=1$,  то из условия (\ref{Aydyad1}) и
(\ref{Aydyad2}) следует, что
\begin{equation}\label{usl-inv}
\displaystyle\int\limits_{\mathbb{R}^{n}}\rho(t;{\bf
x})d\Gamma({\bf x})=\int\limits_{\mathbb{R}^{n}}\rho(0;{\bf
y})d\Gamma({\bf y})=1,
\end{equation}
т. е. для  случайной функции $\rho(t;{\bf x})=\rho(t;{\bf
x};\omega)$ существует случайный функционал, сохраняющий постоянное значение:
\begin{equation}\label{ro1}
\displaystyle\int\limits_{\mathbb{R}^{n}}\rho(t;{\bf
x})d\Gamma({\bf x})=1,
\end{equation}
который можно рассматривать как стохастический объем.

Тогда \eqref{Aydyad1} с условиями \eqref{Aydyad2} и \eqref{Aydyad3} можно рассматривать как стохастический интегральный инвариант.

\begin{dfn} Неотрицательную функцию $\rho(t;{\bf x})$ будем называть
стохастическим ядром или стохастической плотностью стохастического интегрального
инварианта, если выполняются равенства {\rm(\ref{Aydyad1})}, \eqref{Aydyad2} и \eqref{Aydyad3}.
\end{dfn}

\begin{remark}
Понятие интегрального инварианта для системы обыкновенных дифференциальных уравнений было известно ранее, например, оно рассматривалось В.\,И.~Зубовым в {\rm\cite[\S 8
]{Zubov_82}}. Однако существенное отличие, позволившее в {\rm\cite{D_89,D_02}} рассматривать инвариантность случайного объема на основе ядра интегрального инварианта, состоит в том, что в \eqref{Aydyad1} присутствует функциональный множитель. Таким образом, понятие ядра интегрального инварианта для системы обыкновенных дифференциальных уравнений можно рассматривать как частный случай введенного в {\rm\cite{D_89,D_02}}, если взять $f(t;{\bf x})=1$ и, исключив случайность, заданную винеровскими и пуассоновскими процессами, рассматривать интегрирование по детерминированному объему.
\end{remark}

Определим соотношения, при которых функция $\rho(t;{\bf x})$ для
произвольной дважды дифференцируемой функции $f(t;{\bf x})$ будет
ядром интегрального инварианта.

Для случайной функции $f(t;{\bf x}(t))$, где ${\bf x}(t)$ --
решение уравнения (\ref{Ayd01.1}), запишем обобщенную формулу Ито
{\rm{\cite[с.\, 271-272]{GS_68}}}:
\begin{equation}\label{Ayd2}
\begin{array}{c}
 \displaystyle d_{t}f(t;{\bf x}(t))=\Bigl[ \frac{\partial f(t;{\bf x}(t)) }{\partial t}+
 \sum\limits_{i=1}^{n}a_{i}(t;{\bf x}(t))\frac{\partial f(t;{\bf x}(t)) }{\partial x_{i}}+\Bigr.\\
\Bigl.+ \displaystyle  \frac{1}{2} \sum\limits_{i=1}^{n}
 \sum\limits_{j=1}^{n}
 \sum\limits_{k=1}^{m}b_{i\,k}(t;{\bf x}(t))b_{j\,k}(t;{\bf x}(t))
 \frac{\partial^{\,2} f(t;{\bf x}(t)) }{\partial x_{i}x_{j}}\Bigr]dt +\\
 +\displaystyle \sum\limits_{i=1}^{n}
 \sum\limits_{k=1}^{m}b_{i\,k}(t;{\bf x}(t))\frac{\partial f(t;{\bf x}(t)) }{\partial x_{i}} dw_{k}(t)+\\
 +\displaystyle \int
 \Bigl[ f\left(t;{\bf x}(t)+g(t;{\bf x}(t);\gamma)\right)- f(t;{\bf x}(t))
 \Bigr]\nu(dt;d\gamma).
 \end{array}
\end{equation}

\begin{remark}
В дальнейшем, на протяжении главы, для упрощения записей, будем иметь в
виду, что по индексам, встречающимся дважды, проводится
суммирование (без использования знака суммы).
\end{remark}

Продифференцируем по $t$ обе части равенства (\ref{Aydyad1}),
учитывая, что в левой части $f(t;{\bf x})$  -- детерминированная
функция, а в правой $f(t;{\bf x}(t;{\bf y}))$ -- случайный
процесс. Имеем:
\begin{equation*}
\begin{array}{c}
  \displaystyle\int\limits_{\mathbb{R}^{n}}\Bigl(f(t;{\bf
x})d_{t}\rho(t;{\bf x})+\rho(t;{\bf x})\frac{\partial f(t;{\bf
x})}{\partial t}dt\Bigr)d\Gamma({\bf
x})= \int\limits_{\mathbb{R}^{n}}\rho(0;{\bf y})d_{t}f(t;{\bf x}(t;{\bf
y}))d\Gamma({\bf y}).
\end{array}
\end{equation*}
Тогда, в силу (\ref{Aydyad1}) и (\ref{Ayd2}), получаем:
\begin{equation}\label{Aydyad5}
\begin{array}{c}
  \displaystyle\int\limits_{\mathbb{R}^{n}}\Bigl(f(t;{\bf
x})d_{t}\rho(t;{\bf x})+\rho(t;{\bf x})\frac{\partial f(t;{\bf
x})}{\partial t}dt\Bigr)d\Gamma({\bf
x})=\\
=\displaystyle\int\limits_{\mathbb{R}^{n}}\rho(0;{\bf y})d_{t}
f(t;{\bf x}(t;{\bf y}))d\Gamma({\bf y})
=\displaystyle\int\limits_{\mathbb{R}^{n}}\rho(t;{\bf
x})d_{t}f(t;{\bf
x})d\Gamma({\bf x})= \\
=\displaystyle\int\limits_{\mathbb{R}^{n}}d\Gamma({\bf
x})\rho(t;{\bf x})\cdot \biggl[ \Bigl[\frac{\partial f(t;{\bf
x}(t)) }{\partial t}+
 a_{i}(t;{\bf x})\frac{\partial f(t;{\bf x}) }{\partial x_{i}}+\Bigr.\biggr.\\
\Bigl.+ \displaystyle  \frac{1}{2} b_{i\,k}(t;{\bf
x})b_{j\,k}(t;{\bf x})
 \frac{\partial^{\,2} f(t;{\bf x}) }{\partial x_{i} \partial x_{j}}\Bigr]dt
 +\displaystyle
 b_{i\,k}(t;{\bf x})\frac{\partial f(t;{\bf x}) }{\partial x_{i}} dw_{k}(t)+\\
\biggl. +\displaystyle \int
 \Bigl[ f\left(t;{\bf x}+g(t;{\bf x};\gamma)\right)- f(t;{\bf x})
 \Bigr]\nu(dt;d\gamma)\biggr].
\end{array}
\end{equation}
Рассмотрим интеграл
\begin{equation}\label{Aydy1}
I=\displaystyle\int\limits_{\mathbb{R}^{n}}d\Gamma({\bf
x})\rho(t;{\bf x})
 f\left(t;{\bf x}+g(t;{\bf x};\gamma)\right).
\end{equation}
Сделаем замену переменных:
\begin{equation}\label{Aydx1y}
{\bf x}+g(t;{\bf x};\gamma)={\bf y}.
\end{equation}

\begin{remark}
 Область интегрирования в \eqref{Aydy1} разбиваем на подобласти, которые не пересекаются, и на которых реализуется однозначное соответствие в \eqref{Aydy1} между ${\bf x}$ и ${\bf y}$. В дальнейшем используем единое обозначение решения ${\bf x}^{-1}(t;{\bf y};\gamma)$  равенства \eqref{Aydx1y}, относительно ${\bf x}$  для этих подмножеств.
\end{remark}

Тогда, учитывая переход от ${\bf x}$ к ${\bf y}$, получаем:
$
{\bf x}={\bf y}-g(t;{\bf x};\gamma)={\bf y}-g(t;{\bf
x}^{-1}(t;{\bf y};\gamma);\gamma).
$
Следовательно, интеграл $I$ примет вид:
\begin{equation}\label{Aydy2}
\begin{array}{c}
  I=\displaystyle \int\limits_{\mathbb{R}^{n}}d\Gamma({\bf
y})\rho\left(t;{\bf y}-g(t;{\bf x}^{-1}(t;{\bf
y};\gamma);\gamma)\right)
  \cdot f(t;{\bf y})\cdot D\left( {\bf x}^{-1}(t;{\bf y};\gamma)  \right),
\end{array}
\end{equation}
где $D\left( {\bf x}^{-1}(t;{\bf y};\gamma)\right) $ -- якобиан
перехода от ${\bf x}$ к ${\bf y}$.

С учетом (\ref{Aydy2}) и интегрирования по всему пространству
$\mathbb{R}^{n}$, что дает возможность формальной замены
обозначения переменной интегрирования, получаем:
\begin{equation}\label{Aydy3}
\begin{array}{c}
  \displaystyle\int\limits_{\mathbb{R}^{n}}d\Gamma({\bf
x})\rho(t;{\bf x})
 \displaystyle \int
 \Bigl[ f\left(t;{\bf x}+g(t;{\bf x};\gamma)\right)- f(t;{\bf x})
 \Bigr]\nu(dt;d\gamma)= \\
  =\displaystyle\int\limits_{\mathbb{R}^{n}}d\Gamma({\bf
x})\rho(t;{\bf x})
 \displaystyle \int
f\left(t;{\bf x}+g(t;{\bf
x};\gamma)\right)\nu(dt;d\gamma)-
\\
-\displaystyle\int\limits_{\mathbb{R}^{n}}d\Gamma({\bf
x})\rho(t;{\bf x})
 \displaystyle \int  f(t;{\bf
x})
\nu(dt;d\gamma)=
\end{array}
\end{equation}
$$
\begin{array}{c}
=\displaystyle \int\limits_{\mathbb{R}^{n}}d\Gamma({\bf x})
  \int   \Bigl(\Bigr.\rho\left(t;{\bf x}-g(t;{\bf x}^{-1}(t;{\bf
x};\gamma);\gamma)\right) \cdot
 f(t;{\bf x})\cdot\\
  \cdot D\left( {\bf x}^{-1}(t;{\bf x};\gamma)  \right)\Bigl.\Bigr) \nu(dt;d\gamma)-\\
 -\displaystyle\int\limits_{\mathbb{R}^{n}}d\Gamma({\bf
x})\rho(t;{\bf x})
 \displaystyle \int  f(t;{\bf x}) \nu(dt;d\gamma)=\\
=\displaystyle \int\limits_{\mathbb{R}^{n}}d\Gamma({\bf
x})f(t;{\bf x})\int
\Bigl[\rho\left(t;{\bf x}-g(t;{\bf x}^{-1}(t;{\bf
x};\gamma);\gamma)\right)
 \cdot \Bigr.\\
 \cdot\Bigl. D\left( {\bf x}^{-1}(t;{\bf x};\gamma)  \right)-\rho(t;{\bf
 x})\Bigr]\nu(dt;d\gamma)
\end{array}
$$

Учитывая (\ref{Aydyad2}), вычислим следующие интегралы, используя
интегрирование по частям:
\begin{equation}\label{Aydi01}
\displaystyle\int\limits_{\mathbb{R}^{n}}d\Gamma({\bf
x})\rho(t;{\bf x})a_{i}(t;{\bf x})\frac{\partial f(t;{\bf x})
}{\partial x_{i}},
\end{equation}
\begin{equation}\label{Aydi02}
\displaystyle\int\limits_{\mathbb{R}^{n}}d\Gamma({\bf
x})\rho(t;{\bf x})b_{i\,k}(t;{\bf x})\frac{\partial f(t;{\bf x})
}{\partial x_{i}},
\end{equation}
\begin{equation}\label{Aydi03}
\displaystyle\int\limits_{\mathbb{R}^{n}}d\Gamma({\bf
x})\rho(t;{\bf x})b_{i\,k}(t;{\bf x})b_{j\,k}(t;{\bf x})
 \frac{\partial^{\,2} f(t;{\bf x}) }{\partial x_{i} \partial x_{j}}.
\end{equation}
Для (\ref{Aydi01}) имеем:
\begin{equation*}
\begin{array}{c}
  \displaystyle\int\limits_{\mathbb{R}^{n}}d\Gamma({\bf
x})\rho(t;{\bf x})a_{i}(t;{\bf x})\frac{\partial f(t;{\bf x})
}{\partial x_{i}}
=\displaystyle\int\limits_{-\infty}^{+\infty}\prod\limits_{j=1}^{n}dx_{j}\rho(t;{\bf
x})a_{i}(t;{\bf x})\frac{\partial f(t;{\bf x}) }{\partial x_{i}}= \\
  =\displaystyle\int\limits_{-\infty}^{+\infty}\prod\limits_{j=1}^{n-1}dx_{j}
  \int\limits_{-\infty}^{+\infty}\rho(t;{\bf
x})a_{i}(t;{\bf x})\frac{\partial f(t;{\bf x}) }{\partial
x_{i}}dx_{i}.
\end{array}
\end{equation*}
С учетом (\ref{Aydyad3}) вычислим внутренний интеграл с помощью интегрирования по частям:
\begin{equation*}
\begin{array}{c}
\displaystyle\int\limits_{-\infty}^{+\infty}\rho(t;{\bf
x})a_{i}(t;{\bf x})\frac{\partial f(t;{\bf x}) }{\partial
x_{i}}dx_{i}=\\
=\rho(t;{\bf x})a_{i}(t;{\bf x})f(t;{\bf
x})\biggl|_{-\infty}^{+\infty}-\displaystyle\int\limits_{-\infty}^{+\infty}f(t;{\bf
x})\displaystyle\frac{\partial \left( \rho(t;{\bf x})a_{i}(t;{\bf
x})\right) }{\partial
x_{i}}dx_{i} =
\end{array}
\end{equation*}
$$
\begin{array}{c}
=-\displaystyle\int\limits_{-\infty}^{+\infty}f(t;{\bf
x})\displaystyle\frac{\partial \left( \rho(t;{\bf x})a_{i}(t;{\bf
x})\right) }{\partial x_{i}}dx_{i} .
\end{array}
$$
Следовательно,
\begin{equation}\label{Aydi10}
\begin{array}{c}
\displaystyle\int\limits_{\mathbb{R}^{n}}d\Gamma({\bf
x})\rho(t;{\bf x})a_{i}(t;{\bf x})\frac{\partial f(t;{\bf x})
}{\partial x_{i}}
=\displaystyle\int\limits_{-\infty}^{+\infty}\prod\limits_{j=1}^{n}dx_{j}\rho(t;{\bf
x})a_{i}(t;{\bf x})\frac{\partial f(t;{\bf x}) }{\partial x_{i}}=
\end{array}
\end{equation}
$$
  =\displaystyle\int\limits_{-\infty}^{+\infty}\prod\limits_{j=1}^{n-1}dx_{j}
  \int\limits_{-\infty}^{+\infty}\rho(t;{\bf
x})a_{i}(t;{\bf x})\frac{\partial f(t;{\bf x}) }{\partial
x_{i}}dx_{i}=
$$
$$
=-\displaystyle\int\limits_{-\infty}^{+\infty}\prod\limits_{j=1}^{n}dx_{j}
\displaystyle\int\limits_{-\infty}^{+\infty}f(t;{\bf
x})\displaystyle\frac{\partial \left( \rho(t;{\bf x})a_{i}(t;{\bf
x})\right) }{\partial
x_{i}}dx_{i} =
$$
$$
=
-\displaystyle\int\limits_{\mathbb{R}^{n}}d\Gamma({\bf
x})f(t;{\bf x})\displaystyle\frac{\partial \left( \rho(t;{\bf
x})a_{i}(t;{\bf x})\right) }{\partial x_{i}}.
$$
Аналогичным образом вычислим второй интеграл (\ref{Aydi02})  и,
применяя дважды интегрирование по частям, вычислим третий интеграл
(\ref{Aydi03}):
\begin{equation}\label{Aydi20}
\begin{array}{c}
  \displaystyle\int\limits_{\mathbb{R}^{n}}d\Gamma({\bf
x})\rho(t;{\bf x})b_{i\,k}(t;{\bf x})\frac{\partial f(t;{\bf x})
}{\partial x_{i}}
  =-\displaystyle\int\limits_{\mathbb{R}^{n}}d\Gamma({\bf
x})f(t;{\bf x})\frac{\partial\rho(t;{\bf x})b_{i\,k}(t;{\bf x})
}{\partial x_{i}},
\end{array}
\end{equation}

\begin{equation}\label{Aydi30}
\begin{array}{c}
  \displaystyle\int\limits_{\mathbb{R}^{n}}d\Gamma({\bf
x})\rho(t;{\bf x})b_{i\,k}(t;{\bf x})b_{j\,k}(t;{\bf x})
 \frac{\partial^{\,2} f(t;{\bf x}) }{\partial x_{i} \partial x_{j}}=\\
 =\displaystyle\int\limits_{\mathbb{R}^{n}}d\Gamma({\bf x})f(t;{\bf
x})
 \frac{\partial^{\,2} \left(\rho(t;{\bf x})b_{i\,k}(t;{\bf x})b_{j\,k}(t;{\bf x})\right) }{\partial x_{i} \partial
 x_{j}}.
\end{array}
\end{equation}
В (\ref{Aydyad5}) перенесем все в правую часть и с учетом
(\ref{Aydy3}),
 (\ref{Aydi10}), (\ref{Aydi20}) и (\ref{Aydi30}) получим:
\begin{equation}\label{Aydi}
\begin{array}{c}
 0= \displaystyle\int\limits_{\mathbb{R}^{n}}d\Gamma({\bf x})f(t;{\bf
x})\cdot \biggl[-d_{t}\left(\rho(t;{\bf x})\right)-\frac{\partial
f(t;{\bf x}(t)) }{\partial t}+\biggr.\\
+ \displaystyle\frac{\partial f(t;{\bf x}(t)) }{\partial
t}-\frac{\partial\rho(t;{\bf x})b_{i\,k}(t;{\bf x}) }{\partial
x_{i}}dw_{k}(t)+
\end{array}
\end{equation}
$$
+\Bigl(\Bigr.-\displaystyle\frac{\partial \left( \rho(t;{\bf
x})a_{i}(t;{\bf x}(t))\right) }{\partial x_{i}}
+\displaystyle\frac{1}{2}\frac{\partial^{\,2} \left(\rho(t;{\bf
x})b_{i\,k}(t;{\bf x})b_{j\,k}(t;{\bf x})\right) }{\partial x_{i} \partial x_{j}}\Bigl.\Bigr)dt +
$$
$$
 \biggl.+\displaystyle\int  \Bigl[\rho\left(t;{\bf
x}-g(t;{\bf x}^{-1}(t;{\bf x};\gamma);\gamma)\right)
  \cdot D\left( {\bf x}^{-1}(t;{\bf x};\gamma)  \right)-\rho(t;{\bf
 x})\Bigr]\nu(dt;d\gamma)\biggr].
$$
Чтобы равенство (\ref{Aydi}) имело место для любой локально
ограниченной функции $f(t;{\bf x})$,  имеющей ограниченные
производные второго порядка, интегральный инвариант $\rho(t;{\bf
x})$ должен являться решением стохастического уравнения
\begin{equation}\label{Aydii}
\begin{array}{c}
 d_{t}\rho(t;{\bf x})=
-\displaystyle\frac{\partial\rho(t;{\bf x})b_{i\,k}(t;{\bf x})
}{\partial x_{i}}dw_{k}(t) +\Bigl(-\displaystyle\frac{\partial
\left( \rho(t;{\bf x})a_{i}(t;{\bf x})\right) }{\partial
x_{i}}+
\end{array}
\end{equation}
$$
+\displaystyle\frac{1}{2}\frac{\partial^{\,2} \left(\rho(t;{\bf
x})b_{i\,k}(t;{\bf x})b_{j\,k}(t;{\bf x})\right) }{\partial x_{i} \partial x_{j}}\Bigr)dt +
$$
$$
 +\displaystyle\int  \Bigl[\rho\left(t;{\bf
x}-g(t;{\bf x}^{-1}(t;{\bf x};\gamma);\gamma)\right)
 \cdot D\left( {\bf x}^{-1}(t;{\bf x};\gamma)  \right)-\rho(t;{\bf
 x})\Bigr]\nu(dt;d\gamma).
$$
При этом должны выполняться условия
\begin{equation}\label{Aydii1}
\begin{array}{c}
\rho(t;{\bf x})\Bigl|_{t=0}=\rho(0;{\bf x})\in C_{0}^{2},  \\
\lim\limits_{|{\bf x}|\to \infty}\rho(0;{\bf x})=0, \ \ \ \
\displaystyle\lim\limits_{|{\bf x}|\to
\infty}\frac{\partial\rho(0;{\bf x})}{\partial x_{i}}=0.
 \end{array}
\end{equation}
Таким образом, получены условия инвариантности стохастического объема и доказана следующая
\begin{theorem}\label{thro1}
Пусть ${\bf x}(t)$,
${\bf x}\in\mathbb{R}^{n}$, решение системы обобщенных стохастических
дифференциальных уравнений Ито
\begin{equation*}
\begin{array}{c}
  dx_{i}(t)= \displaystyle a_{i}\Bigl(t;{\bf x}(t) \Bigr)\, dt
  +\sum\limits_{k=1}^{m}
b_{i,k}\Bigl(t;{\bf x}(t)\Bigr)\, dw_{k}(t) +\\
+
\displaystyle\int g_{i}\Bigl(t;{\bf x}(t);\gamma\Bigr)\nu(dt;d\gamma), \\
  {\bf x}(t)={\bf x}\Bigl(t;{\bf x}(0) \Bigr) \Bigr|_{t=0}={\bf
  x}(0), \ \ \ \ \ \ i=\overline{1,n}, \ \ \ \ \ \ t\geq 0,
\end{array}
\end{equation*}
коэффициенты которого удовлетворяют условию $\mathfrak{L})$;
${\bf w} (t)$ -- $m$-мерный винеровский процесс,
$\nu(t;\Delta\gamma)$ -- однородная по $t$ нецентрированная мера
Пуассона и $\rho(t;{\bf x})$ -- случайная функция, измеримая
относительно потока $\sigma$-алгебр
$\Bigl\{\mathcal{F}\Bigr\}_{0}^{T}$,
$\mathcal{F}_{t_{1}}\subset\mathcal{F}_{t_{2}}$, $t_{1}<t_{2}$,
согласованного с процессами ${\bf w} (t)$ и $\nu(t;\Delta\gamma)$
и  относительно любой функции $f(t;{\bf x})$ из класса $\mathfrak{S}$ локально ограниченных функций, имеющих ограниченные вторые
производные по ${\bf x}$. Функция $\rho(t;{\bf x})$ является стохастическом ядром стохастического интегрального инварианта
\eqref{Aydyad1} для
произвольной локально ограниченной функции $f(t;{\bf x})\in \mathfrak{S}$, если она является решением системы \eqref{Aydii} обобщенных СДУ Ито, удовлетворяющим начальным условиям \eqref{Aydii1}.
\end{theorem}

\subsection*{2.2. Обобщенная формула Ито -- Вентцеля }
\addcontentsline{toc}{subsection}{2.2. Обобщенная формула Ито -- Вентцеля}

Равенство (\ref{Aydyad1}) с условиями (\ref{Aydyad2}), (\ref{Aydyad3}) и уравнение (\ref{Aydii}) для ядра
интегрального инварианта соответствовали случаю детерминированной
функции $f(t;{\bf x})$.  Положим выполнение аналогичного равенства  для случайной функции $F(t;{\bf x})$:
\begin{equation}\label{Aydz3p}
\begin{array}{c}
  \displaystyle\int\limits_{\mathbb{R}^{n}}\rho(t;{\bf x})F(t;{\bf
x})d\Gamma({\bf x})=
  \displaystyle\int\limits_{\mathbb{R}^{n}}\rho(0;{\bf y})F(t;{\bf x}(t;{\bf y}))d\Gamma({\bf y}),
\end{array}
\end{equation}
где ${\bf x}(t;{\bf y})$ -- решение системы СДУ (\ref{Ayd01.1}).

Рассмотрим сложный случайный процесс $F\left(t;{\bf x}(t;{\bf
y})\right)\in \mathbb{R}^{n_{o}}$, где ${\bf x}(t;{\bf y})$ --
решение системы СДУ (\ref{Ayd01.1}), а процесс $F(t;{\bf
x})$ -- решение системы обобщенных СДУ Ито:
\begin{equation}\label{Aydz1}
\begin{array}{c}
  d_{t}F(t;{\bf x})=Q(t;{\bf x})dt+D_{k}(t;{\bf x})dw_{k}(t)
+
  \displaystyle\int \nu(dt;d\gamma)G(t;{\bf x};\gamma).
  \end{array}
\end{equation}
Относительно случайных функций $Q(t;{\bf x})$, $D_{k}(t;{\bf
x})$, $G(t;{\bf x};\gamma)$, определенных на пространстве
$\mathbb{R}^{n_{o}}$, предполагаем, что они непрерывны и
ограничены вместе со своими производными по всем переменным,
измеримые относительно потока $\sigma$-алгебр $\mathcal{F}_{t}$,
согласованного с процессами ${\bf w}(t)$ и $\nu(t;\Delta\gamma)$
из (\ref{Ayd01.1}).

Опираясь на уравнение для стохастического интегрального
инварианта, построим правило дифференцирования
для сложного случайного процесса $F\left(t;{\bf x}(t;{\bf
y})\right)$ \cite{11_KchOboz}.

Рассмотрим интеграл
$
\displaystyle\int\limits_{\mathbb{R}^{n}}\rho(t;{\bf x})F(t;{\bf x})d\Gamma({\bf x})
$.
Поскольку интегрирование проводится по пространству
$\mathbb{R}^{n}$, на котором задан процесс ${\bf x}(t)$,
используем равенство (\ref{Aydz3p}), записав его в виде (поменяв части равенства местами):
\begin{equation}\label{Aydz3}
\begin{array}{c}
  \displaystyle\int\limits_{\mathbb{R}^{n}}\rho(0;{\bf y})F(t;{\bf x}(t;{\bf y}))d\Gamma({\bf y})=
  \int\limits_{\mathbb{R}^{n}}\rho(t;{\bf x})F(t;{\bf
x})d\Gamma({\bf x}).
\end{array}
\end{equation}
Продифференцируем обе части (\ref{Aydz3}) по $t$:
\begin{equation}\label{Aydz4}
\begin{array}{c}
  \displaystyle\int\limits_{\mathbb{R}^{n}}\rho(0;{\bf y})d_{t}F(t;{\bf x}(t;{\bf y}))d\Gamma({\bf
y})= \\
 \displaystyle\int\limits_{\mathbb{R}^{n}}\Bigl(\Bigr.\rho(t;{\bf x})d_{t}F(t;{\bf x})+F(t;{\bf x})d_{t}\rho(t;{\bf x}) -
D_{k}(t;{\bf x}) \displaystyle\frac{\partial\rho(t;{\bf
x})b_{i\,k}(t;{\bf x}) }{\partial x_{i}}dt\Bigl.\Bigr)d\Gamma({\bf
x}) ,
\end{array}
\end{equation}
Учитывая, что интегрирование идет по пространству $\mathbb{R}^{n}$
после введения замены переменной, то при интегрировании по кривой
${\mathbb R}(\gamma)$ в этом пространстве нужно учитывать произведенную
замену (\ref{Aydx1y}).

Поскольку функция $\rho(t;{\bf x})$ -- ядро интегрального инварианта \eqref{Aydz3p}, применим теорему~\ref{thro1}.  Подставим (\ref{Aydii}) и (\ref{Aydz1}) в правую часть
(\ref{Aydz4}):
$$
I_{1}= \displaystyle \int\limits_{\mathbb{R}^{n}}\Bigl(\rho(t;{\bf
x})d_{t}F(t;{\bf x})+F(t;{\bf x})d_{t}\rho(t;{\bf x})
-D_{k}(t;{\bf x}) \displaystyle\frac{\partial\rho(t;{\bf
x})b_{i\,k}(t;{\bf x})
}{\partial x_{i}}dt\Bigr)d\Gamma({\bf x})=
$$
$$  = \displaystyle\int\limits_{\mathbb{R}^{n}}d\Gamma({\bf x})\biggl\{\rho(t;{\bf x})
  \Bigl[\Bigr.Q(t;{\bf x})dt+D_{k}(t;{\bf x})dw_{k}(t)
+
$$
$$
+\displaystyle\int G(t;{\bf x}+g(t;{\bf
x}^{-1}(t;{\bf
x};\gamma);\gamma)\nu(dt;d\gamma)\Bigl.\Bigr]+
$$
$$
+F(t;{\bf x})\biggl[-\displaystyle\frac{\partial\rho(t;{\bf
x})b_{i\,k}(t;{\bf x}) }{\partial
x_{i}}dw_{k}(t)+\Bigr.
$$
$$ +\displaystyle\Bigl(\Bigr.-\displaystyle\frac{\partial
\left( \rho(t;{\bf x})a_{i}(t;{\bf x})\right) }{\partial
x_{i}}+
\displaystyle\frac{1}{2}\frac{\partial^{\,2} \left(\rho(t;{\bf
x})b_{i\,k}(t;{\bf x})b_{j\,k}(t;{\bf x})\right) }{\partial x_{i} \partial x_{j}}\Bigl.\Bigr)dt +
$$
$$
 +\displaystyle\int  \Bigl[\Bigr.\rho\left(t;{\bf
x}-g(t;{\bf x}^{-1}(t;{\bf x};\gamma);\gamma)\right)
 \cdot  D\left( {\bf x}^{-1}(t;{\bf x};\gamma)  \right)-\rho(t;{\bf
 x})\Bigl.\Bigr]\nu(dt;d\gamma)\Bigl.\Bigl.\biggr]-
 $$
 $$
 -D_{k}(t;{\bf
x}) \displaystyle\frac{\partial\rho(t;{\bf x})b_{i\,k}(t;{\bf x})
}{\partial x_{i}}dt\biggr\}.
$$
Приведем подобные:
\begin{equation}\label{Aydz5}
\begin{array}{c}
I_{1}=\displaystyle\int\limits_{\mathbb{R}^{n}}d\Gamma({\bf
x})\biggl(\rho(t;{\bf x})Q(t;{\bf x}) -F(t;{\bf
x})\displaystyle\frac{\partial \left( \rho(t;{\bf x})a_{i}(t;{\bf
x})\right) }{\partial x_{i}}+\biggr.\\
-D_{k}(t;{\bf x}) \displaystyle\frac{\partial\rho(t;{\bf
x})b_{i\,k}(t;{\bf x}) }{\partial x_{i}}
+\displaystyle\frac{1}{2}F(t;{\bf x})\frac{\partial^{\,2}
\left(\rho(t;{\bf x})b_{i\,k}(t;{\bf x})b_{j\,k}(t;{\bf x})\right)
}{\partial x_{i}
\partial x_{j}}\biggl.\biggr)dt+\\
+\displaystyle\int\limits_{\mathbb{R}^{n}}d\Gamma({\bf
x})\biggl(\rho(t;{\bf x})D_{k}(t;{\bf x}) -F(t;{\bf
x})\displaystyle\frac{\partial \left( \rho(t;{\bf
x})b_{i\,k}(t;{\bf x})\right) }{\partial
x_{i}}\biggr)dw_{k}(t)+\\
+\displaystyle\int\limits_{\mathbb{R}^{n}}d\Gamma({\bf
x})\biggl[\rho(t;{\bf x}) \cdot\int G(t;{\bf
x}+g(t;{\bf x}^{-1}(t;{\bf
x};\gamma);\gamma))\nu(dt;d\gamma) +\Bigr.\\
+F(t;{\bf
x})\displaystyle\int \Bigl[\rho\left(t;{\bf
x}-g(t;{\bf x}^{-1}(t;{\bf x};\gamma);\gamma)\right)
 \cdot D\left( {\bf x}^{-1}(t;{\bf x};\gamma)  \right)-\rho(t;{\bf
 x})\Bigr]\nu(dt;d\gamma)\Bigl.\biggr].
\end{array}
\end{equation}
В силу (\ref{Aydi10}), (\ref{Aydi20}), (\ref{Aydi30}) имеем:
\begin{equation*}
\begin{array}{c}
\displaystyle\int\limits_{\mathbb{R}^{n}}d\Gamma({\bf x})F(t;{\bf x})\frac{\partial \left(\rho(t;{\bf x})a_{i}(t;{\bf
x})\right) }{\partial x_{i}} =
-\displaystyle\int\limits_{\mathbb{R}^{n}}d\Gamma({\bf
x})\rho(t;{\bf x})a_{i}(t;{\bf x})\displaystyle\frac{\partial
 F(t;{\bf x})
 }{\partial x_{i}},
\end{array}
\end{equation*}
\begin{equation*}
\begin{array}{c}
  \displaystyle\int\limits_{\mathbb{R}^{n}}d\Gamma({\bf
x})D_{k}(t;{\bf x})\frac{\partial\left(\rho(t;{\bf
x})b_{i\,k}(t;{\bf x})\right) }{\partial x_{i}}=
  -\displaystyle\int\limits_{\mathbb{R}^{n}}d\Gamma({\bf
x})\rho(t;{\bf x})b_{i\,k}(t;{\bf x})\frac{\partial D_{k}(t;{\bf
x}) }{\partial x_{i}},
\end{array}
\end{equation*}
\begin{equation*}
\begin{array}{c}
\displaystyle\int\limits_{\mathbb{R}^{n}}d\Gamma({\bf x})F(t;{\bf x})
 \frac{\partial^{\,2} \left(\rho(t;{\bf x})b_{i\,k}(t;{\bf x})b_{j\,k}(t;{\bf x})\right) }{\partial x_{i} \partial
 x_{j}}=
  \\
 = \displaystyle\int\limits_{\mathbb{R}^{n}}d\Gamma({\bf
x})\rho(t;{\bf x})b_{i\,k}(t;{\bf x})b_{j\,k}(t;{\bf x})
 \frac{\partial^{\,2} F(t;{\bf x}) }{\partial x_{i} \partial
 x_{j}}.
\end{array}
\end{equation*}
Вычислим последний интеграл в сумме (\ref{Aydz5}):
\begin{equation}\label{Ayde6}
\begin{array}{c}
  I_{2}=\displaystyle\int\limits_{\mathbb{R}^{n}}d\Gamma({\bf
x})F(t;{\bf x})
  \cdot\rho\Bigl(t;{\bf x}-g(t;{\bf x}^{-1}(t;{\bf
x};\gamma);\gamma)\Bigr)D({\bf x}^{-1}(t;{\bf x};\gamma)).
\end{array}
\end{equation}
Введем замену переменных:
$$
\begin{array}{c}
  {\bf x}-g(t;{\bf x}^{-1}(t;{\bf x};\gamma);\gamma)={\bf y}, \\
  {\bf x}={\bf y}+g(t;{\bf x}^{-1}(t;{\bf x};\gamma);\gamma)={\bf y}+g(t;{\bf
  y};\gamma).
\end{array}
$$
Обозначим якобиан перехода от вектора ${\bf x}$ к вектору ${\bf
y}$ через \linebreak $D_{o}({\bf x}^{-1}(t;{\bf y};\gamma))$. Тогда, в силу
(\ref{Aydy1}) и дальнейшей формальной замены обозначения
переменной интегрирования, получаем:
\begin{equation}\label{Ayde7}
\begin{array}{c}
  I_{2}=\displaystyle\int\limits_{\mathbb{R}^{n}}d\Gamma({\bf y})F\Bigl(t;{\bf y}+g(t;{\bf
  y};\gamma)\Bigr)
  \cdot\rho(t;{\bf y})D_{o}({\bf x}^{-1}(t;{\bf y};\gamma))D({\bf
x}^{-1}(t;{\bf x};\gamma))= \\
 =\displaystyle\int\limits_{\mathbb{R}^{n}}d\Gamma({\bf y})F\Bigl(t;{\bf y}+g(t;{\bf
  y};\gamma)\Bigr)\rho(t;{\bf y})
  =\displaystyle\int\limits_{\mathbb{R}^{n}}d\Gamma({\bf x})F\Bigl(t;{\bf x}+g(t;{\bf
  x};\gamma)\Bigr)\rho(t;{\bf x}).
\end{array}
\end{equation}
В результате, правая часть (\ref{Aydz4}) примет вид:
\begin{equation}\label{Ayde8}
\begin{array}{c}
I_{1}=\displaystyle\int\limits_{\mathbb{R}^{n}}d\Gamma({\bf
x})\rho(t;{\bf x})\biggl\{\Bigl( D_{k}(t;{\bf x}) +b_{i\,k}(t;{\bf
x})\frac{\partial F(t;{\bf x}) }{\partial x_{i}}
\Bigr)dw_{k}(t)+\biggr.\\+
 \Bigl(Q(t;{\bf x})+a_{i}(t;{\bf
x})\displaystyle\frac{\partial
 F(t;{\bf x})
 }{\partial x_{i}}+ b_{i\,k}(t;{\bf
x})\frac{\partial
 D_{k}(t;{\bf x})
 }{\partial x_{i}}\Bigr.
 +\\
 +\Bigl.\displaystyle\frac{1}{2}b_{i\,k}(t;{\bf x})b_{j\,k}(t;{\bf x})
 \frac{\partial^{\,2} F(t;{\bf x}) }{\partial x_{i} \partial
 x_{j}}\Bigr)dt+\\
 +\displaystyle\int G(t;{\bf x}+g(t;{\bf
x}^{-1}(t;{\bf x};\gamma);\gamma)\nu(dt;d\gamma)+\\
+
 \displaystyle\int \Bigl[F\Bigl(t;{\bf x}+g(t;{\bf
  x};\gamma)\Bigr)
  -F(t;{\bf
 x})\Bigr]\nu(dt;d\gamma)\biggl.\biggr\}.
 \end{array}
\end{equation}
В (\ref{Aydz4}) перенесем все в правую часть и с учетом
(\ref{Aydyad1}) получаем:
\begin{equation*}
\begin{array}{c}
0=\displaystyle\int\limits_{\mathbb{R}^{n}}d\Gamma({\bf
y})\rho(0;{\bf y})\biggl\{-d_{t}F(t;{\bf x}(t;{\bf
y}))+\biggr.\\
+\Bigl( D_{k}(t;{\bf x}(t;{\bf y}))+b_{i\,k}(t;{\bf x}(t;{\bf
y}))\displaystyle\frac{\partial F(t;{\bf
x}) }{\partial x_{i}} \Bigr)dw_{k}(t)+\\
+ \Bigl(Q(t;{\bf x}(t;{\bf y}))+a_{i}(t;{\bf
x})\displaystyle\frac{\partial
 F(t;{\bf x})
 }{\partial x_{i}}
 +b_{i\,k}(t;{\bf
x})\frac{\partial
 D_{k}(t;{\bf x})
 }{\partial x_{i}}+\Bigr.\\
 +\Bigl.\displaystyle\frac{1}{2}b_{i\,k}(t;{\bf x}(t;{\bf y}))b_{j\,k}(t;{\bf x}(t;{\bf y}))
 \frac{\partial^{\,2} F(t;{\bf x}(t;{\bf y})) }{\partial x_{i} \partial
 x_{j}}\Bigr)dt+\\
 +\displaystyle\int G(t;{\bf x}(t;{\bf y})+g(t;{\bf
  x}(t;{\bf y});\gamma)\nu(dt;d\gamma)+\\
 +
 \displaystyle\int \Bigl[F\Bigl(t;{\bf x}(t;{\bf y})+g(t;{\bf
  x}(t;{\bf y});\gamma)\Bigr)-
F(t;{\bf
 x}(t;{\bf y}))\Bigr]\nu(dt;d\gamma)\biggl.\biggr\}.
 \end{array}
\end{equation*}
Таким образом, получаем обобщенную формулу Ито -- Вентцеля  -- формулу \eqref{GrindEQ__2_5_2_} для нецентрированной пуассоновской меры:
\begin{equation*}
\begin{array}{c}
d_{t}F(t;{\bf x}(t;{\bf y}))=\displaystyle\Bigl(
D_{k}(t;{\bf x}(t;{\bf y}))
+b_{i\,k}(t;{\bf x}(t;{\bf
y}))\frac{\partial F(t;{\bf
x}(t;{\bf y})) }{\partial x_{i}} \Bigr)dw_{k}(t)+\\
 +\Bigl(Q(t;{\bf x}(t;{\bf y}))+a_{i}(t;{\bf
x}(t;{\bf y}))\displaystyle\frac{\partial
 F(t;{\bf x}(t;{\bf y}))
 }{\partial x_{i}}+
  \end{array}
  \end{equation*}
$$
 +b_{i\,k}(t;{\bf
x}(t;{\bf y}))\displaystyle\frac{\partial
 D_{k}(t;{\bf x}(t;{\bf y}))
 }{\partial x_{i}}+\Bigr.
 $$
 $$
 +\Bigl.\displaystyle\frac{1}{2}b_{i\,k}(t;{\bf x}(t;{\bf y}))b_{j\,k}(t;{\bf x}(t;{\bf y}))
 \frac{\partial^{\,2} F(t;{\bf x}(t;{\bf y})) }{\partial x_{i} \partial
 x_{j}}\Bigr)dt+
$$
\begin{equation*}
\begin{array}{c}
 +\displaystyle\int G(t;{\bf x}(t;{\bf y})+g(t;{\bf
  x}(t;{\bf y});\gamma)\nu(dt;d\gamma)+\\
 +
 \displaystyle\int \Bigl[F\Bigl(t;{\bf x}(t;{\bf y})+g(t;{\bf
  x}(t;{\bf y});\gamma)\Bigr)
  -F(t;{\bf
 x}(t;{\bf y}))\Bigr]\nu(dt;d\gamma).
 \end{array}
\end{equation*}

Представленное B. Oksendal  и T. Zhang в \cite{Oks_07} (2007) обобщения формулы Ито -- Вентцеля произведено для уравнения, содержащего только  пуасссоновскую составляющую (без винеровской) и только для скалярного случайного процесса, поэтому его можно рассматривать как частный случай полученной формулы \eqref{GrindEQ__2_5_2_}.

\begin{remark}В {\rm\cite[c.\,24]{D_02}} было введено понятие стохастического первого
интеграла для центрированной пуассоновской меры и полученные
условия для его существования учитывают необходимость задания
плот\-но\-сти интенсивности пуассоновского распределения  в отличие от
предложенного в данной работе. Таким образом, безразлично, каков
вероятностный закон имеют интенсивности пуассоновских скачков. Это
обстоятельство является очень важным для дальнейших применений, в
частности, для построения программных управлений {\rm\cite{11_KchUpr}}.
\end{remark}

Предложенное обобщение формулы Ито -- Вентцеля  и понятия
стохастического первого интеграла {\rm{\cite{D_02}}} позволяет,
как отмечено в {\rm{\cite{07_ChOboz}}}, строить программные
управления стохастических динамических систем, подверженных
случайным возмущениям, вызванным винеровскими возмущениями и
пуассоновскими скачками
{\rm{\cite{11_KchUpr}}}.

\renewcommand{\refname}{Библиографический список}


\begin{thebibliography}{99}
\addcontentsline{toc}{section}{Библиографический список \dotfill}





\bibitem{D_84}
Дубко В. А. Интегральные инварианты для одного класса систем
стохастических дифференциальных уравнений // Докл. АН УССР, Сер.~А, №1, 1984. -- С. 18--21.

\bibitem{D_89}
Дубко~В.\,А. Вопросы теории и применения
стохастических дифференциальных урав\-не\-ний.~--  Владивосток:
ДВНЦ АН СССР, 1989.~---  185~c.

\bibitem{D_02}
Дубко~В.\,А. Открытые эволюционирующие системы.// "<Вiдкритi еволюцiонуючi системи">,
мiжнар. наук.-практ. конф. (2002, Київ) Перша мiжнародна
науково-практична конференцiя "<Вiдкритi еволюцiонуючi системи">,
(26-27 квiт. 2002 р.) (Додаток). К., ВНЗ ВМУРоЛ, 2002.~---
С.~14-31.


\bibitem{GS_68}
Гихман И.И.,
Скороход  А.В. Введение в теорию случайных процессов. --  М. :
Наука, 1965. -- 654 с.

\bibitem{Gelfand_00}
Гельфанд~И.М., Шилов Г.Е. Обобщенные функции и действия над ними. - М.: Добросвет, 2000. - 412 с.


\bibitem{Zubov_82}
Зубов~В.\,И. Динамика управляемых систем: Учебное пособие для вузов. M. : Высшая школа, 1982. -- 285 c.

\bibitem{11_KchUpr}
Карачанская~Е.\,В. Построение программных управлений
с вероятностью~1 для динамической системы с пуассоновскими
возмущениями.
// Вестник Тихоокеанского госуниверситета, No 2 (21) 2011. -- \linebreak С.~51--60.

\bibitem{11_KchOboz}
{Карачанская Е.В. Об одном обобщении формулы Ито -- Вентцеля. // Обозрение прикладной и промышленной математики,
т.18, вып.~2, 2011.~-- С.~494--496}

\bibitem{Rozovsky_73}
{Розовский Б.Л. О формуле Ито -- Вентцеля. // Вестник МГУ. Сер. матем. механ. -- 1973. -- № 1. --  С.~26--32.}

\bibitem{07_ChOboz}
{Чалых Е.В. Программное управление с вероятностью 1 для открытых
систем. // Обозрение прикладной и промышленной математики, т.14,
вып.~2, 2007.~-- С.~253-254.}

\bibitem{Oks_07}
{Oksendal, B. and Zhang, T. The Ito-Ventcel formula and forward
stochastic differential equation driven by Poisson random
measures.
//  Osaka J. Math. 44 (2007), pp. 207--230.}


\end{thebibliography}
\end{document}